\documentclass[english,12pt,reqno]{smfart}
\usepackage{amssymb,amsmath,amstext,amsthm,amsfonts}
\usepackage{graphicx}
\usepackage[ansinew]{inputenc}
\usepackage{a4wide}
\newcommand{\length}{\operatorname{length}}
\newcommand{\leb}{\operatorname{Leb}}

\newcommand{\dist}{\operatorname{dist}}

\newcommand{\lip}{\operatorname{Lip}}

\begin{document}

\newcommand{\mcup}{\mbox{$\bigcup$}}
\newcommand{\mcap}{\mbox{$\bigcap$}}

\def \RR {{\mathbb R}}
\def \ZZ {{\mathbb Z}}
\def \NN {{\mathbb N}}
\def \PP {{\mathbb P}}
\def \TT {{\mathbb T}}
\def \II {{\mathbb I}}
\def \JJ {{\mathbb J}}

\def \vare {\varepsilon }

 \def \cf {\mathcal{F}}
 \def \cm {\mathcal{M}}
 \def \cn {\mathcal{N}}
 \def \cq {\mathcal{Q}}
 \def \cp {\mathcal{P}}
 \def \cb {\mathcal{B}}
 \def \cc {\mathcal{C}}
 \def \cs {\mathcal{S}}
 \def \bc {\mathcal{B}}
 \def \hc {\mathcal{C}}

\newcommand{\dem}{\begin{proof}}
\newcommand{\cqd}{\end{proof}}

\newcommand{\qand}{\quad\text{and}\quad}

\newtheorem{theorem}{Theorem}
\newtheorem{corollary}{Corollary}

\newtheorem*{Maintheorem}{Main Theorem}

\newtheorem{maintheorem}{Theorem}
\renewcommand{\themaintheorem}{\Alph{maintheorem}}
\newcommand{\cmt}{\begin{maintheorem}}
\newcommand{\fmt}{\end{maintheorem}}

\newtheorem{maincorollary}[maintheorem]{Corollary}
\renewcommand{\themaintheorem}{\Alph{maintheorem}}
\newcommand{\cmc}{\begin{maincorollary}}
\newcommand{\fmc}{\end{maincorollary}}

\newtheorem{T}{Theorem}[section]
\newcommand{\ct}{\begin{T}}
\newcommand{\ft}{\end{T}}

\newtheorem{Corollary}[T]{Corollary}
\newcommand{\cco}{\begin{Corollary}}
\newcommand{\fco}{\end{Corollary}}

\newtheorem{Proposition}[T]{Proposition}
\newcommand{\cpr}{\begin{Proposition}}
\newcommand{\fpr}{\end{Proposition}}

\newtheorem{Lemma}[T]{Lemma}
\newcommand{\cle}{\begin{Lemma}}
\newcommand{\fle}{\end{Lemma}}

\newtheorem{Sublemma}[T]{Sublemma}
\newcommand{\csle}{\begin{Lemma}}
\newcommand{\fsle}{\end{Lemma}}

\theoremstyle{definition}

\newtheorem{Remark}[T]{Remark}
\newcommand{\cre}{\begin{Remark}}
\newcommand{\fre}{\end{Remark}}

\newtheorem{Definition}[T]{Definition}
\newcommand{\cd}{\begin{Definition}}
\newcommand{\fd}{\end{Definition}}

\title[Gibbs-Markov structures for partially hyperbolic attractors]{Gibbs-Markov structures and limit laws for partially hyperbolic attractors with mostly expanding central direction}

\author{Jos\'e F. Alves}
\address{Departamento de Matem\'atica Pura, Faculdade de Ci\^encias da Universidade do Porto\\
Rua do Campo Alegre 687, 4169-007 Porto, Portugal}
\email{jfalves@fc.up.pt}
\urladdr{http://www.fc.up.pt/cmup/home/jfalves}

\author{Vilton  Pinheiro}
\address{Departamento de Matem\'atica, Universidade Federal da Bahia\\
Av. Ademar de Barros s/n, 40170-110 Salvador, Brazil.}
\email{viltonj@ufba.br}

\date{\today}

\thanks{Work carried out at the Federal University of
Bahia and University of Porto. JFA was partially supported by
FCT through CMUP and by POCI/MAT/61237/2004. VP was partially
supported by PROCAD/Capes and by POCI/MAT/61237/2004}



\begin{abstract}
We consider a partially hyperbolic set $K$ on a Riemannian
manifold $M$ whose tangent space  splits as $T_K
M=E^{cu}\oplus E^{s}$, for which the centre-unstable direction $E^{cu}$ expands non-uniformly
on some local unstable disk.
We show that under these assumptions $f$ induces a Gibbs-Markov structure. Moreover, the decay of
the return time function
 can be controlled in terms of the time typical points need to achieve
 some uniform expanding behavior in the centre-unstable direction.
 As an application of the main result we obtain certain
 rates for decay of correlations, large deviations, an almost sure invariance principle and the
 validity of the Central Limit Theorem.
\end{abstract}

\maketitle

\setcounter{tocdepth}{2}

\tableofcontents


\section{Introduction}

Remarkable advances in the study of  dynamical systems, specially the statistical properties of those with chaotic behavior,  
have been  achieved through the idea of \emph{inducing}. Roughly speaking, this consists of replacing the initial dynamical system by another one whose dynamical features are easier to understand and from which one can recover much information on the initial system. This idea goes back to the 70's where Markov partitions have been used to study the statistical properties of uniformly hyperbolic dynamical systems via conjugations to shifts. Since then, a main goal in dynamical systems theory  is to enlarge that strategy to wider classes of systems.

A main achievement in this direction has been attained by Young in \cite{Y1,Y2}. In these works she developed an abstract framework, the so-called \emph{Young towers}, that proved usefulness  in a systematic  treatment of several classes of dynamical systems, including Axiom~A attractors, piecewise hyperbolic maps, billiards with convex scatterers, logistic maps, intermittent maps and Hénon-type attractors. The latter case has actually been treated by Benedicks and Young in~\cite{BY2}. A preponderant role in this context is played by  \emph{Gibbs-Markov structures}, which may be understood as a generalization of the classical Markov partitions and are naturally associated to an inducing scheme that gives rise to a Young tower. Many statistical properties on the dynamics of these induced structures can be recovered from the \emph{Gibbs-Markov map} that one obtains quotienting out by stable leaves. Gibbs-Markov maps constitute themselves object of great dynamical interest; see e.g. \cite{AD01,AD02} for Aaronson and Denker contributions on the statistical properties of Gibbs-Markov maps, whose main ideas go back to Aaronson's book \cite{Aar} on infinite ergodic theory.

A Gibbs-Markov structure is characterized by a suitable region of the phase space partitioned into  subsets (possibly infinitely many) each of which with a  given \emph{return time}. Comparing to the classical  Markov partitions, a main difference lies on the possibility of infinitely many return times, as long as the measure of points with larger and larger returns decays to zero. The flexibility conferred by the chance of arbitrarily large return times is a fundamental step towards  applications to non-uniformly hyperbolic dynamical systems, where large waiting times are needed to reach good expansion rates in the centre-unstable direction of some points. In certain cases we are able to determine the speed at which the return times decay, given in terms of the time that generic points need to achieve the good expansion rates.

Dynamical systems that do not fit the class of uniformly hyperbolic ones but combine non-uniform expanding/contracting central directions with other directions of uniformly hyperbolic behavior give rise to a wider class of partially hyperbolic dynamical  systems. The  case of partially hyperbolic diffeomorphisms for which the tangent bundle over some attracting set splits as $E^{cs}\oplus E^u$, the sum of one invariant sub-bundle with non-uniformly contracting behavior  and another one with uniformly expanding behavior, has been treated in~\cite{Cas}, where some Gibbs-Markov structures were obtained to deduce decay of correlations and the validity of the Central Limit Theorem for the SRB measures which had been obtained in~\cite{BV}; see also \cite{Car,PS}. Limit theorems for partially hyperbolic systems of this type were also obtained in \cite{D}.

Our main goal in this work is to prove the existence of Gibbs-Markov structures in the case  of  partially hyperbolic diffeomorphisms with an attracting set over which the tangent bundle splits as $E^{s}\oplus E^{cu}$, the sum of a sub-bundle $E^{s}$ having uniformly contracting behavior with another one $E^{cu}$ having non-uniform expanding behavior. This kind of attracting set has been previously considered in~\cite{ABV}, where the existence of \emph{SRB measures} was established.
The method used in \cite{ABV} is based on the simple idea of iterating forward Lebesgue measure on some centre-unstable disk and obtaining the absolute continuity with respect to Lebesgue measure on local unstable disks of weak* accumulation points. Due to the simplicity of the method, it gives not much information on the properties of such SRB measures. As a byproduct of the machinery that we develop here, we are also able to deduce the existence of the SRB measure. Our method uses deeper knowledge on the geometrical structure of the attractor, thus enabling us to prove the existence of a Gibbs-Markov structure inside it, leading to an inducing scheme.  As an immediate consequence of our main result, combined with others from  \cite{AP2,Mel,MN2}, we easily deduce some statistical properties of these dynamical systems, namely \emph{decay of correlations}, \emph{central limit theorem}, \emph{large deviations} and a \emph{multidimensional almost sure invariant principle}. Let us point out that the case $E^{s}\oplus E^{cu}$ that we consider here is, for our purposes, considerably more difficult to deal with than the dual case $E^{cs}\oplus E^u$,  since the richest part of the dynamics in the neighborhood of an attracting set occurs in the unstable direction, where in our case the expansion is attained just asymptotically.


The final part of the proof of our main result follows the strategy used in \cite{ALP} for non-uniformly expanding endomorphisms (non-invertible smooth dynamical systems). The argument gives no optimal conclusions outside the polynomial case. The lack of efficiency for exponential or subexponential decays essentially relies on \cite[Proposition 6.1]{ALP} which still has no suitable generalizations. A main achievement in this direction has been obtained in \cite{G} by mean of a different geometrical construction, leading to exponential and subexponential decays of return times. The construction  in  \cite{ALP} can be thought of as being local, in the sense that  the partition is obtained by considering convenient returns of points in a small disk to itself.  The construction performed  in \cite{G} uses instead  a finite global partition of the whole attractor as a starting point. However, this strategy has no natural generalization to the present setting of partially hyperbolic diffeomorphisms, mostly due to the non-compactness of unstable manifolds.

\subsection*{Overview} This paper is organized as follows. In the remaining of this introduction we consider three subsections. In the first one we define the Gibbs-Markov structures. In the second subsection we introduce partial hyperbolicity and state our main result on the existence of Gibbs-Markov structures for certain partially hyperbolic attractors. In the final subsection we define some limit laws and state some statistical consequences of our main result. In Section~\ref{sub.hyptimes} we recall some results on hyperbolic times and bounded distortion from~\cite{ABV}. The construction of Gibbs-Markov structures for partially hyperbolic attractors is performed in Section~\ref{se.Gibbs-Markov structure}. We begin with some results on the recurrence of disks, and then present an algorithmic construction that gives rise to the product structure. Finally, in Subsection~\ref{sec.regularity} we prove some results on the regularity of the stable and unstable foliations. This comprises the generalization of classical results for uniformly hyperbolic attractors to our setting, namely the Hölder continuity of the central-unstable direction and the absolute continuity of the stable foliation. Finally, the estimates on the decay of return times is obtained in Section~\ref{s.metric}.

 \subsection{Gibbs-Markov structures}\label{se.hypstructures}  Here we present the structures which have been introduced in \cite{Y1} and constitute the main object of  our interest.  These structures comprise the dynamical and geometrical essence of the return map to the base of a Young tower.

Let $f:M\to M$ be a diffeomorphism of a Riemannian manifold $M$.  We say that $f$ is $C^{1+}$ if $f$ is $C^1$ and $Df$
is H\"older continuous.
Let $\leb$ denote the
Lebesgue measure on the Borel sets of $M$ associated to the Riemannian structure. Given a submanifold $\gamma\subset M$   we use
$\leb_\gamma$ to denote the Lebesgue measure on $\gamma$ induced by the
restriction of the Riemannian structure to~$\gamma$.

An embedded disk $\gamma\subset M$ is called an {\em unstable
manifold} if $\dist(f^{-n}(x),f^{-n}(y))\to0$ as
$n\to\infty$ for every $x,y\in\gamma$. Similarly, $\gamma$ is called
a {\em stable manifold} if $\dist(f^{n}(x),f^{n}(y))\to0$
 as $n\to\infty$ for every $x,y\in\gamma$.
 Let  $\text{Emb}^1(D^u,M)$ be the space of $C^1$ embeddings from
$D^u$ into $M$. We say that $\Gamma^u=\{\gamma^u\}$ is a
\emph{continuous family of $C^1$ unstable manifolds} if there is a
compact set~$K^s$, a unit disk $D^u$ of some $\RR^n$, and a map
$\Phi^u\colon K^s\times D^u\to M$ such that
\begin{itemize}
\item[$\bullet$]  $\gamma^u=\Phi^u(\{x\}\times D^u)$ is an unstable
manifold;
\item[$\bullet$]  $\Phi^u$ maps $K^s\times D^u$ homeomorphically onto its
image; \item[$\bullet$]  $x\mapsto \Phi^u\vert(\{x\}\times D^u)$ is a
continuous map from $K^s$ to $\text{Emb}^1(D^u,M)$.
\end{itemize}
Continuous families of $C^1$ stable manifolds are defined similarly.
 We say that $\Lambda\subset M$ has a \emph{product
structure} if there exist a continuous family of unstable
manifolds $\Gamma^u=\{\gamma^u\}$ and a continuous family of
stable manifolds $\Gamma^s=\{\gamma^s\}$ such that
\begin{itemize}
    \item[$\bullet$] $\Lambda=(\cup \gamma^u)\cap(\cup\gamma^s)$;
    \item[$\bullet$] $\dim \gamma^u+\dim \gamma^s=\dim M$;
    \item[$\bullet$] each  $\gamma^s$ meets each $\gamma^u$  in exactly one point;
    \item[$\bullet$] stable and unstable manifolds meet transversally with angles bounded
    away from~0.
\end{itemize}

Let $\Lambda\subset M$ have a
product structure whose associated defining families are $\Gamma^s$
and $\Gamma^u$. A subset $\Lambda_0\subset \Lambda$ is called an
{\em $s$-subset} if $\Lambda_0$ also has a hyperbolic product
structure, and its defining families, $\Gamma_0^s$ and $\Gamma_0^u$,
can be chosen with $\Gamma_0^s\subset\Gamma^s$ and
$\Gamma_0^u=\Gamma^u$; {\em $u$-subsets} are defined similarly.
Given $x\in\Lambda$, let $\gamma^{*}(x)$ denote the element of
$\Gamma^{*}$ containing $x$, for $*=s,u$. For each $n\ge 1$ we let
$(f^n)^u$ denote the restriction of the map $f^n$ to
$\gamma^u$-disks and $\det D(f^n)^u$ denote   the Jacobian of
$D(f^n)^u$.

We require that the  product structure  satisfies several properties that we explicit below in (P$_1$)-(P$_5$).
From here on we assume that $C>0$ and $0<\beta<1$ are constants  depending only  on~$f$ and~$\Lambda$.

\begin{enumerate}
    \item[\bf\quad
    (P$\bf_1$)] \emph{Markov}: there are pairwise disjoint $s$-subsets $\Lambda_1,\Lambda_2,\dots\subset\Lambda$ such
    that
    \begin{enumerate}
 \item $\leb_{\gamma}\big((\Lambda\setminus\cup\Lambda_i)\cap\gamma\big)=0$ on each $\gamma\in\Gamma^u$;
 \item for each $i\in\NN$ there is $R_i\in\NN$ such that $f^{R_i}(\Lambda_i)$ is $u$-subset,
         and for all $x\in \Lambda_i$
        \begin{quote}$
        \displaystyle
         f^{R_i}(\gamma^s(x))\subset \gamma^s(f^{R_i}(x))\qand
         f^{R_i}(\gamma^u(x))\supset \gamma^u(f^{R_i}(x)).
         $
         \end{quote}
         This allows us to define the \emph{return time} function $R\colon \Lambda\to \NN$ as $R\vert\Lambda_i=R_i$.
    \end{enumerate}
\end{enumerate}


 \begin{enumerate}
\item[\bf(P$\bf_2$)] \emph{Contraction on  $\Gamma^s$}:
   for all $y\in\gamma^s(x)$ and $ n\ge
 1$
 \begin{quote}$\displaystyle\dist(f^n(y),f^n(x))\le C\beta^n .$
 \end{quote}

 \end{enumerate}


 \begin{enumerate}
\item[\bf(P$\bf_3$)] \emph{Backward contraction on $\Gamma^u$}:  for all $x,y\in\Lambda_i$ with $ y\in\gamma^u(x) $, and $ 0\le n<R_i$
\begin{quote}
$\displaystyle\dist(f^n(y),f^n(x))\le C\beta^{{R_i}-n} \dist(f^{R_i}(x),f^{R_i}(y)).$\end{quote}
 \end{enumerate}


\begin{enumerate}
 \item[\bf(P$\bf_4$)] \emph{Bounded distortion}:
for all $x,y\in\Lambda_i$ with $ y\in\gamma^u(x) $
   \begin{quote}
    $\displaystyle\log\frac{\det D(f^{R_i})^u(x)}{\det D(f^{R_i})^u(y)}\le
    C\dist(f^{R_i}(x),f^{R_i}(y))^{\eta}.$
    \end{quote}
 \end{enumerate}

\begin{enumerate}
\item[\bf(P$\bf_5$)] \emph{Regularity of the foliations}:
\begin{enumerate}
 \item
 for all $y\in\gamma^s(x)$ and $ n\ge 0$
 \begin{quote}$\displaystyle
 \log \prod_{i=n}^\infty\frac{\det Df^u(f^i(x))}{\det Df^u(f^i(y))}\le C\beta^{n};
 $\end{quote}
 \item given $\gamma,\gamma'\in\Gamma^u$, define
$\phi\colon\gamma\cap\Lambda\to\gamma'\cap\Lambda$  as
$\phi(x)=\gamma^s(x)\cap \gamma$. Then $\phi$ is absolutely continuous and
        \begin{quote}$
        \displaystyle \frac{d(\phi_*^{-1}\leb_{\gamma})}{d\leb_{\gamma'}}(x)=
        \prod_{i=0}^\infty\frac{\det Df^u(f^i(x))}{\det
        Df^u(f^i(\phi(x)))}.
        $\end{quote}
\end{enumerate}
\end{enumerate}
See Section~\ref{sec.regularity} for a precise definition of absolute continuity. A set with a product structure for which properties (P$\bf_1$)-(P$\bf_5$) above hold will be called a \emph{Gibbs-Markov structure}.

%
%

\subsection{Partially hyperbolic attractors}
\label{ss.statement}

 Let
$K\subset M$ be a compact {\em invariant} set for a $C^1$ diffeomorphism $f\colon M\to M$, meaning that
$f(K)=K$. We say
that $K$ has a \emph{dominated splitting} if there exists a
continuous $Df$-invariant splitting $T_K M=E^{cs}\oplus
E^{cu}$ and
$0<\lambda<1$ such that for some choice of a Riemannian metric on
$M$ 
 \begin{equation}\label{domination1}
    \|Df \mid E^{cs}_x\|
\cdot \|Df^{-1} \mid E^{cu}_{f(x)}\| \le\lambda,\quad\text{for
all $x\in K$.}
\end{equation}
We call $E^{cs}$ the {\em centre-stable bundle} and $E^{cu}$ the
{\em centre-unstable bundle}.
We say that $K$ is {\em
partially hyperbolic\/} if it has a dominated splitting $T_K
M=E^{cs}\oplus E^{cu}$ for which $E^{cs}$ is \emph{uniformly contracting}:
there is $0<\lambda<1$ such that for some choice of a Riemannian
metric on $M$
\begin{equation}
\label{contraction} \|Df \mid E^{cs}_x\|\le \lambda,\quad\text{for
all $x\in K$.}
\end{equation}
Fixing some small number $c>0$, we say that $f$ is
\emph{non-uniformly expanding in the central-unstable direction} at
a point  $x\in K$ if
\begin{equation}
\label{NUE}\tag{NUE} \limsup_{n\to+\infty} \frac{1}{n}
    \sum_{j=1}^{n} \log \|Df^{-1} \mid E^{cu}_{f^j(x)}\|<-c.
\end{equation}
We
 define the \emph{expansion time} function
\begin{equation}\label{exptime}
    \mathcal E(x) = \min\left\{N\ge 1\colon  \frac{1}{n}
\sum_{i=0}^{n-1} \log \|Df^{-1}\mid E^{cu}_{f^{i}(x)}\| <
-c, \quad \forall\,n\geq N\right\}.
\end{equation}
Note that $\mathcal E(x) $ is finite for the points $x\in K$ satisfying \ref{NUE}.

In this work we consider partially hyperbolic sets of the same type of those considered in~\cite{ABV}, for which
the centre-stable
direction is uniformly contracting and the central-unstable
direction is non-uniformly expanding. To highlight the uniform
contraction in the centre-stable direction we shall write $E^s$
instead of $E^{cs}$.


\cmt \label{t:Markov towers} Let \( f: M\to M \) be a  \( C^{1+} \)
diffeomorphism and let $K\subset M$ be a transitive partially hyperbolic set
such that $T_K M=E^s\oplus E^{cu}$. Assume that for some local
unstable disk $D\subset K$ and \( \tau>0\) one has
    $
   \leb_D\{\mathcal E> n\} = O(n^{-\tau}).
    $
Then there exists $\Lambda\subset K$ with a Gibbs-Markov structure. Moreover,
 $
     \leb_\gamma\{ R> n\} = O(n^{-\tau})
 $
 for any  $\gamma\in \Gamma^u$.
  \fmt

\cre
The existence of a disk~$D$ for which \ref{NUE} is satisfied Lebesgue almost everywhere can actually hold under very mild conditions. Indeed, as shown in
\cite[Theorem~A]{AP}, if $K$ attracts a positive Lebesgue measure set of points for which \ref{NUE} holds, then~$K$ contains some local unstable disk~$D$ for which \ref{NUE} holds for $\leb_D$ almost every point.
\fre

\cre Under the assumptions of Theorem~\ref{t:Markov towers} we are able to say more about the  set $\Lambda$ with the product structure. Actually, our construction gives that the set $\Lambda$ itself coincides with the union of the leaves in $\Gamma^u$. This is not always the case, e. g. for Hénon attractors $\Lambda$ is a Cantor set; see \cite{BY2}.
\fre

An open class of diffeomorphisms for which $K=M$ is partially hyperbolic and satisfies the assumptions of Theorem~\ref{t:Markov towers} can
be found in \cite[Appendix A]{ABV}. The calculations in \cite{ABV} give that for such diffeomorphisms one has $\leb_D\{\mathcal E> n\}$ decaying exponentially fast, which then implies that  $\leb_\gamma\{ R> n\}$ has \emph{super-polynomial decay}, meaning that  it decays faster than any polynomial.  Transitivity of the diffeomorphisms in that class has been proved in~\cite{T}.


\subsection{Limit laws}
An $f$-invariant Borel probability $\mu$ in $M$ is called an \emph{SRB measure} if, for a positive Lebesgue
measure set of points $x\in M$,
\begin{equation}
\lim_{n\to+\infty} \frac{1}{n}
     \sum_{j=0}^{n-1} \varphi(f^j(x)) = \int \varphi\,d\mu,
\quad\text{for any continuous } \varphi:M\to\RR.
\end{equation}
%
Defining the
\emph{correlation function} of observables $\varphi, \psi\colon M\to \RR$ as
\[
\mathcal C_{n}(\varphi, \psi) = \left|\int (\varphi \circ f^{n}) \psi d\mu - \int \varphi d\mu \int \psi
d\mu\right|,
\]
it is sometimes possible to obtain specific rates for which $C_{n}(\varphi, \psi)$ decays to 0 as $n$ tends to infinity, at least for certain classes of observables with some regularity.
Note that taking the observables as characteristic functions of Borel sets we are led to the classical definition of \emph{mixing}.

The next two corollaries follow from Theorem~\ref{t:Markov towers} together with \cite[Theorem~B]{AP2} and \cite[Theorem C]{AP2}; see also 
\cite[Remark~2.4]{AP2}.

\cmc[Decay of Correlations]\label{c:decay} Let \( f: M\to M \) be a  \( C^{1+} \)
diffeomorphism and let $K\subset M$ be a transitive partially hyperbolic set
such that $T_K M=E^s\oplus E^{cu}$. Assume that there is a local
unstable disk $D\subset K$ and \( \tau>1\) such that
    $
   \leb_D\{\mathcal E> n\}=O(n^{-\tau}).
    $
Then some power of \( f \) has an SRB  measure \( \mu \) and  
    $
    \mathcal C_{n}(\varphi,\psi) = O(n^{-\tau+1 })
    $
    for H\"older
continuous  \(\varphi,\psi\colon M\to\RR \).
\fmc

The existence of the SRB measure for $f$ has already been proved in \cite[Theorem A]{ABV}. 
In general, we can only assure that the correlation decay holds for some power of $f$. However, if the return times associated to the elements of the Gibbs-Markov structure given by Theorem~\ref{t:Markov towers} are relatively prime, i.e. gcd$\{R_i\}=1$, then the same conclusion holds with respect to~$f$. \emph{For simplicity, from here on we assume that gcd$\{R_i\}=1$.} Otherwise, the same conclusions hold for some power of~$f$. 


In the next result we obtain conditions for the validity of the
\emph{Central Limit
Theorem (CLT)},  which states that the deviation of
the average values of an observable along an orbit from the
asymptotic average is given by a Normal Distribution:
given any Hölder continuous $\varphi\colon M\to\RR$  which is not a
\emph{coboundary} (\( \varphi\neq \psi\circ f - \psi \) for any \( \psi \in L^2\))
there exists \( \sigma>0 \) such that 
\[
\frac{1}{\sqrt n}\sum_{j=0}^{n-1}\left(\varphi\circ f^{j}-\int\varphi d\mu
\right) \stackrel{\text{distr}}\longrightarrow N(0,\sigma), \quad\text{as $n\to\infty$}.
\]

 \cmc[Central Limit Theorem]\label{c:tlc} Let \( f: M\to M \) be a  \( C^{1+} \)
diffeomorphism and let $K\subset M$ be a transitive partially hyperbolic set
with $T_K M=E^s\oplus E^{cu}$. Assume that there is a local
unstable disk $D\subset K$ and \( \tau>2\) such that
    $
   \leb_D\{\mathcal E> n\}=O(n^{-\tau}).
    $
    Then CLT holds for any Hölder continuous $\varphi\colon M\to\RR$  which is not a
coboundary.
 \fmc

Given a  Hölder continuous observable $\varphi\colon M\to \RR$ and $\epsilon>0$, we define the \emph{large deviation} of the time average with respect to the mean of $\varphi$ as
 $$\mathcal D_n (\varphi,\epsilon)=\mu\left(\left\{x\in M\colon \left|\frac{1}{ n}\sum_{j=0}^{n-1}\varphi(f^{j}(x))-\int\varphi d\mu
\right|>\epsilon\right\}\right).
$$
The next result is a consequence of Theorem~\ref{t:Markov towers} together with \cite[Theorem 4.2]{MN1} and \cite[Theorem 1.2]{Mel}.

 \cmc[Large Deviations]\label{c:ld} Let \( f: M\to M \) be a  \( C^{1+} \)
diffeomorphism and let $K\subset M$ be a transitive partially hyperbolic set
with $T_K M=E^s\oplus E^{cu}$. Assume that there is a local
unstable disk $D\subset K$ and \( \tau>2\) such that
    $
   \leb_D\{\mathcal E> n\}= \mathcal O(n^{- \tau}).
    $
Then, for any Hölder continuous $\varphi\colon M\to\RR$ and any $\epsilon > 0$, we have
 $
 \mathcal D_n (\varphi,\epsilon)=O(n^{-\tau+1}) .
 $
 \fmc

Given $d\ge 1$ and a Hölder continuous $\varphi\colon M\to \RR^d$, we denote $\varphi_n=\sum_{j=0}^{n-1}(\varphi\circ f^j-\int\varphi\, d\mu)$, for each $n\ge 1$.
We say that the sequence $\{\varphi_n\}_n$ satisfies a \emph{$d$-dimensional
almost sure invariance principle} (\emph{ASIP})
if there exists $\lambda > 0$ and  a probability
space supporting a sequence of random variables $\{\varphi^*_n\}_n$ and a $d$-dimensional Brownian motion
$W(t)$ such that
 \begin{enumerate}
\item  $\{\varphi_n\}_n $ and $ \{\varphi^*_n\}_n $ are equally distributed;
\item  $\varphi^*_n = W(n) + O(n^{1/2-\lambda})$, as $n\to\infty$, almost everywhere.
 \end{enumerate}
 The ASIP is said to be \emph{nondegenerate} if the Brownian motion $W(t)$ has nonsingular
covariance matrix $\Sigma$. For the dynamical systems considered in this paper, there is a closed subspace
$Z$ of infinite codimension in the space of all (piecewise) H\"older $\varphi\colon M\to \RR^d$ such
that $\Sigma$ is nonsingular whenever $\varphi\notin Z$; see \cite[Remark 1.2]{MN2} and \cite[Section 4.3]{HM}.

The next result follows from Theorem~\ref{t:Markov towers} and \cite[Theorem 1.6]{MN2}, and it generalizes Corollary~\ref{c:tlc} above.

 \cmc[Almost Sure Invariance Principle]\label{c:asip} Let \( f: M\to M \) be a  \( C^{1+} \)
diffeomorphism and let $K\subset M$ be a transitive partially hyperbolic set
with $T_K M=E^s\oplus E^{cu}$. Assume that there is a local
unstable disk $D\subset K$ and \( \tau>2\) such that
    $
   \leb_D\{\mathcal E> n\}= O(n^{-\tau}).
    $
Then $\{\varphi_n\}_n$ satisfies an ASIP for any Hölder continuous $\varphi\notin Z$.
 \fmc


%
%

\section{Preliminaries}\label{sub.hyptimes}

With the only exception of Lemma~\ref{le:grow1}, the material contained in this section comes from~\cite{ABV}.
Our aim is to recall some results on the H\"older continuity of the tangent
direction of centre-unstable submanifolds, to introduce hyperbolic times and recall their main properties.

 We fix continuous extensions of the
two subbundles $E^{s}$ and $E^{cu}$ to some neighborhood $U$ of
$K$ that we denote by $\tilde{E}^{s}$ and $\tilde{E}^{cu}$. We
do not require these extensions to be invariant under $Df$. Given
$0<a<1$, we define the {\em centre-unstable cone field
$C_a^{cu}=\left(C_a^{cu}(x)\right)_{x\in U}$ of width $a$\/} by
\begin{equation}
\label{e.cucone} C_a^{cu}(x)=\big\{v_1+v_2 \in
\tilde{E}_x^{s}\oplus \tilde{E}_x^{cu} \text{\ such\ that\ }
\|v_1\| \le a \|v_2\|\big\}.
\end{equation}
We define the {\em centre-stable cone field
$C_a^{s}=\left(C_a^{s}(x)\right)_{x\in U}$ of width $a$\/} in a
similar way, simply reversing the roles of the subbundles in
(\ref{e.cucone}). We fix $a>0$ and $U$ small enough so that, up to
slightly increasing $\lambda<1$, the domination condition
\eqref{domination1}  remains valid for any pair of vectors in the
two cone fields:
$$
\|Df(x)v^{s}\|\cdot\|Df^{-1}(f(x))v^{cu}\|
\le\lambda\|v^{s}\|\,\|v^{cu}\|
$$
for every $v^{s}\in C_a^{s}(x)$, $v^{cu}\in C_a^{cu}(f(x))$, and
any point $x\in U\cap f^{-1}(U)$. Note that the centre-unstable
cone field is positively invariant:
 $$Df(x) C_a^{cu}(x)\subset
C_a^{cu}(f(x)),\quad\text{whenever $x,f(x)\in U$.}$$ Actually, the
domination property together with the invariance of
$E^{cu}=\tilde{E}^{cu} \vert K$ imply that
$$
Df(x) C_a^{cu}(x) \subset C_{\lambda a}^{cu}(f(x))
                  \subset C_a^{cu}(f(x)),
$$
for every $x\in K$, and this extends to any $x\in U\cap
f^{-1}(U)$ just by continuity.



We say that an embedded $C^1$ submanifold $N\subset U$ is {\em
tangent to the centre-unstable cone field\/}, if the tangent
subspace to $N$ at each point $x\in N$ is contained in the
corresponding cone $C_a^{cu}(x)$. Then $f(N)$ is also tangent to
the centre-unstable cone field, if it is contained in $U$, by the
domination property. The tangent bundle of $N$ is said to be {\em
H\"older continuous} if $x \mapsto T_x N$ defines a H\"older
continuous section from $N$ to the corresponding Grassman bundle
of $M$. Given a $C^1$ submanifold $N\subset U$, we define
\begin{equation}
\label{e.kappa} \kappa(N)=\inf\{C>0:\text{the tangent bundle of $N$
is $(C,\zeta)$-H\"older}\}.
\end{equation}
The next result contains all the information we need on the
H\"older control of the tangent direction and its proof can be
found in \cite[Corollary 2.4]{ABV}.

\cpr \label{c.curvature} There exists $C_1>0$ such that, given any
$C^1$ submanifold $N\subset U$ tangent to the centre-unstable cone
field, then
\begin{enumerate}
  \item there exists $n_0\ge 1$ such that $\kappa(f^n(N)) \le
C_1$ for every $n\ge n_0$ such that $f^k(N) \subset U$ for all $0\le
k \le n$;
  \item  if $\kappa(N) \le C_1$, then $\kappa(f^n(N)) \le
C_1$ for every $n\ge 1$ such that $f^k(N)\subset U$ for all $0\le k
\le n$;
\item if $N$ and $n$ are as in the previous item, then the functions
$$
J_k: f^k(N)\ni x \mapsto \log |\det \big(Df \mid T_x f^k(N)\big)|,
\quad\text{$0\le k \le n$},
$$
are $(L_1,\zeta)$-H\"older continuous with $L_1>0$ depending only on
$C_1$ and $f$.
\end{enumerate}
   \fpr

The notion we introduce next allows us to derive uniform expansion and
bounded distortion estimates from the non-uniform expansion assumption in the
centre-unstable direction.

\cd \label{d.hyperbolic1} Given $0<\sigma<1$, we say that $n$ is a
{\em $\sigma$-hyperbolic time\/} for  $x\in K$ if
$$
\prod_{j=n-k+1}^{n}\|Df^{-1} \mid E^{cu}_{f^{j}(x)}\| \le \sigma^k,
\qquad\text{for all $1\le k \le n$.}
$$
\fd
 In particular, if $n$ is a $\sigma$-hyperbolic time for $x$,
then $Df^{-k} \mid E^{cu}_{f^{n}(x)}$ is a contraction for every
$1\le k \le n$:
\begin{equation}\label{contra}
    \|Df^{-k} \mid E^{cu}_{f^{n}(x)}\| \le \prod_{j=n-k+1}^{n}\|Df^{-1}
\mid E^{cu}_{f^{j}(x)}\| \le \sigma^{k}.
\end{equation}

If $a>0$ is sufficiently small 
and we choose $\delta_1>0$ also small so that the
$\delta_1$-neighborhood of $K$ is contained in $U$, then, by
continuity,
\begin{equation}
\label{e.delta1} \|Df^{-1}(f(y)) v \| \le \frac{1}{\sqrt\sigma}
\|Df^{-1}|E^{cu}_{f(x)}\|\,\|v\|,
\end{equation}
whenever $x\in K$, $\dist(x,y)\le \delta_1$, and $v\in
C^{cu}_a(y)$.

Given any disk $\Delta\subset M$, we use $\dist_\Delta(x,y)$ to
denote the distance between $x,y\in \Delta$, measured along
$\Delta$. The distance from a point $x\in \Delta$ to the boundary
of $\Delta$ is $\dist_\Delta(x,\partial \Delta)=
\inf_{y\in\partial \Delta}\dist_\Delta(x,y)$. The
next result has essentially been proved in \cite[Lemma 2.7]{ABV}; see \cite[Lemma 4.2]{AP} for a detailed proof.

\cle \label{l.contraction} Let $0<\delta<\delta_1$ and $\Delta\subset U$ be a $C^1$ disk of
radius $\delta$ tangent to the
centre-unstable cone field.
Then, there is $n_0\ge1$ such that for $x\in \Delta\cap K$ with
$\dist_\Delta(x,\partial \Delta)\ge \delta/2$ and $n \ge n_0$ a
$\sigma$-hyperbolic time for $x$ there is a neighborhood $V_n$ of
$x$ in $\Delta$ such that:
\begin{enumerate}
    \item $f^{n}$ maps $V_n$ diffeomorphically onto a
centre-unstable disk  of radius $\delta_1$ around  $f^{n}(x)$;
    \item for every $1\le k
\le n$ and $y, z\in V_n$, $$
\dist_{f^{n-k}(V_n)}(f^{n-k}(y),f^{n-k}(z)) \le
\sigma^{k/2}\dist_{f^n(V_n)}(f^{n}(y),f^{n}(z)).$$
\end{enumerate}
\fle

We call the sets \( V_n \)  \emph{hyperbolic pre-balls} and their
images \( f^{n}(V_n) \)  \emph{hyperbolic balls}. Notice that the
latter are indeed centre-unstable balls of radius \( \delta_1 \).
The next result follows from Proposition~\ref{c.curvature} and Lemma~\ref{l.contraction} above exactly as in the proof of \cite[Proposition 2.8]{ABV}.

\cco \label{p.distortion} There exists $C_2>1$ such that given a
disk $\Delta$ as in Lemma~\ref{l.contraction} with $\kappa(\Delta)
\le C_1$, and given any hyperbolic pre-ball $V_n\subset \Delta$ with
$n\ge n_0$, then for all $y,z\in V_n$
$$
\log \frac{|\det Df^{n} \mid T_y \Delta|}
                     {|\det Df^{n} \mid T_z \Delta|}
            \le C_2 \dist_{f^n(D)}(f^{n}(y),f^{n}(y))^\zeta.
$$
 \fco

The next result states the existence of $\sigma$-hyperbolic times
for points satisfying \ref{NUE} and its proof can be found in
\cite[Lemma 3.1, Corollary 3.2]{ABV}.

\cpr\label{l:hyperbolic2}
    There are \( \theta>0 \) and $\sigma>0$ 
     such that for  every $x
    \in  K$ with $\mathcal E(x)\le n$ there exist $\sigma$-hyperbolic times
    $1 \le n_1 < \cdots < n_l \le n$ for \( x \) with $l\ge\theta
    n$.
\fpr

Given $n\ge 1$, we define
 $$
 H_n=\{x\in K\colon \text{ $n$ is a $\sigma$-hyperbolic time for
 $x$ }\}.
 $$
The next result follows from the proposition above as in \cite[Proposition~3.5]{ABV} or \cite[Corollary~2.3]{ALP}. It plays important role in the metric estimates
of Section~\ref{s.metric}.

\begin{Corollary}\label{c:hyperbolic3}
Let $D$ be a local unstable disk for which  \ref{NUE} holds
$\leb_D$ almost everywhere. Given $n\ge1$  and $A\subset
D\setminus\{\mathcal E>n\}$ with $\leb_D(A)>0$ we have
    \[
   \frac{1}{n}\sum_{j=1}^{n}\frac{\leb_D(A\cap H_{j})}{\leb_D(A)} \geq \theta.
   \]
\end{Corollary}

We finish this section with a technical lemma which will be useful in the proof of Proposition~\ref{p.construction}. We take $\delta_s>0$ sufficiently small so that local stable manifolds $W^s_{\delta_s}(x)$ are defined for all points $x\in K$ and
 \begin{equation}\label{eq.c4}
\left|\,\log \|Df^{-1} \mid E^{cu}_{x}\|-\log \|Df^{-1} \mid \tilde
E^{cu}_{y}\| \right|<\frac{c}4,
\end{equation} for all $x\in K$
and $y\in W^s_{\delta_s}(x)$, where $c>0$ is given by \ref{NUE}.

\begin{Lemma}\label{le:grow1}
Given $0<\vare<\delta_1$, there exists  \(  N_{\varepsilon}> 0 \)
such that for every $n\ge n_0$ and every $x\in D\cap H_n$, the ball
of radius $\vare$ centered at $f^n(z)$ inside the hyperbolic ball
$f^n(V_n(x))$ contains a hyperbolic pre-ball $V_k(z)$ with \( k\leq
N_{\varepsilon} \).
\end{Lemma}

\dem Given $n\ge n_0$ and $x\in D\cap H_n$, let $B^n_x=f^n(V_n(x))$ be the
hyperbolic ball associated to $x$ with hyperbolic time $n$. Recall
that $B^n_x$ is a centre-unstable ball of radius
$\delta_1$ around $f^n(x)$. We define the cylinder
  $$C^n_x=\bigcup_{y\in B^n_x}W_{\delta_s}^s(y).$$
Since \ref{NUE} holds for $\leb_D$ almost every point and it
remains valid by forward iteration,
 it follows that $\leb_{B^n_x}$ almost every point in $B^n_x$ also
 satisfies \ref{NUE}.
Given  $z\in C^n_x$, we define
 $$
 \tilde{\mathcal E}(z) = \min\left\{N\ge 1\colon  \frac{1}{n}
\sum_{i=0}^{n-1} \log \|Df^{-1}\mid \tilde E^{cu}_{f^{i}(w)}\| <
-\frac{3c}4, \quad \forall n\geq N \right\}.
 $$
The fact that \ref{NUE} holds for $\leb_{B^n_x}$ almost every
point in $B^n_x$ together with \eqref{eq.c4} imply that
$\tilde{\mathcal E}(z)$ is well defined for $\leb$ almost every
point $z\in C^n_x$. Hence $\tilde{\mathcal E}(z)$ is well defined
for $\leb$ almost every point $z$ belonging to the set
 $$
 C=\bigcup_{n\ge n_0\atop x\in D\cap H_n} C^n_x .
 $$
Using  Lemma~\ref{l.contraction} we may choose \(  n_{\varepsilon}
\) large enough so that any hyperbolic pre-ball of
$\sigma$-hyperbolic time \( n\geq  n_{\varepsilon} \) will have
diameter not exceeding $\vare/2$. Let now $B_x^n(\vare/2)$ denote
the ball of radius $\vare/2$ around $f^n(x)$ inside $B_x^n$, and
take
 $$
 v_\vare=\min_{n\ge n_0 \atop x\in D\cap H_n}
 \big\{\leb\left( \cup_{y\in
 B_x^n(\vare/2)}W_{\delta_s}^s(y)\right)\big\}.
 $$
Since the sizes and tangent directions of hyperbolic balls and
stable disks are uniformly controlled, this minimum $v_\vare$ must
be strictly positive. Hence, as
    \[
     \leb\left\{z\in C\colon \tilde{\mathcal E}(z)>n\right\} \to
     0,
    \quad\text{ when } n\to\infty,
    \]
it is possible to choose $N_\vare\in\NN$ large enough so that
    \begin{equation}\label{isto}
      \leb\left\{z\in C\colon \tilde{\mathcal E}(z)>N_\vare\right\}
    \leq v_\vare.
    \end{equation}
We take $N_\vare$ also satisfying $\theta N_\vare>n_\vare$. By the
choice of $v_\vare$, given $n\ge n_0$ and $x\in D\cap H_n$ there
must be some $z\in W^s_{\delta_s}(y)$ with $y\in B_x^n(\vare/2)$
such that $\tilde{\mathcal E}(z)\le N_\vare$. Using \eqref{eq.c4} we
easily deduce that ${\mathcal E}(z)\le N_\vare$; recall
\eqref{exptime}. Hence, by Proposition~\ref{l:hyperbolic2} there exists
some $\sigma$-hyperbolic time $h$ for $y$ with $\theta N_\vare<h\le
N_\vare$. Since we have taken $\theta N_\vare>n_\vare$, then by
the choice of $n_\vare$ we are done.
\end{proof}

\section{Gibbs-Markov structure}\label{se.Gibbs-Markov structure}

In this section we describe the geometric construction of the product structure. This will be made in three steps. In the first one we prove the existence of a centre-unstable disk (a reference disk) for which forward iterates come back to a neighborhood of itself, and whose projection along stable leaves cover the disk completely. In the next step we use these returns to define a partition on the reference disk. This part of the construction follows ideas from \cite[Section 3]{ALP}. Finally we use the partition on the reference leaf and the returns to define the product structure.

\subsection{Returning disks} Let $D$ be a local unstable disk as in
Theorem~\ref{t:Markov towers}. Diminishing $\delta_1>0$, if
necessary, we may assume that $D$ has radius $\delta_1$. Take
$0<\delta_s<\delta_1/2$ such that points in $K$ have local stable manifolds of radius $\delta_s$. In particular,
these local stable manifolds are contained in $U$; recall
\eqref{e.delta1}.

\cle\label{l.N0q} There are $N_0\ge 1$ and $q\in K$ such
that:
 \begin{enumerate}
 \item $W^s_{\delta_s/2}(q)$ intersects $D$ in a point $p$ with
$\dist_D(p ,\partial D)>\delta_1/2$;
 \item for each centre-unstable disk $\gamma_1^u$ of radius $\delta_1$ centered at a
point in $ K$ there are $0\le j\le N_0$ and a disk
$\gamma_2^u\subset \gamma_1^u$ of radius $\delta_1/2$ centered at a
point $z\in W^s_{\delta_s/4}(f^{-j}(q))$.
 \end{enumerate}
  \fle

\dem We start by observing that there is a constant $\alpha=\alpha(\rho) >0$
with $\alpha\to0$ as $\rho\to0$
for which the following holds:
{\em given $x\in  K$, $\rho>0$  and $y\in  K$ with
$\dist(x,y)<\rho$ having a local unstable disk of radius
$\delta_1$ centered at $y$, then $W^s_{\delta_s}(x)$ intersects
$W^u_{\delta_1}(y)$ in a point $z$ with
$$\dist_{W^s_{\delta_s}(x)}(z,x)<\alpha  \qand\dist_{W^u_{\delta_1}(y)}(z,y)<\delta_1/2.$$}
In particular, such
a point $z$ has a neighborhood of radius $\delta_1/2$ inside
$W^u_{\delta_1}(y)$.

Take $\rho>0$ small so that $4\alpha<\delta_s $.
Since we are assuming $f\vert  K$ transitive, we may fix
$q\in  K$  and
$N_0\in\NN$ such that both $i)$ $W^s_{\delta_s/2}(q)$ intersects
$D$ in a point $p$ with $\dist_D(p ,\partial D)>\delta_1/2$, and
$ii)$  $\{ f^{-N_0}(q),\dots, f^{-1}(q),q\}$ is $\rho$-dense in
$ K$. By $\rho$-dense we mean that any other point in
$ K$ has one of the points in the set above at a distance less than
$\rho$.

Hence, given any centre-unstable disk $\gamma_1^u$ of radius
$\delta_1$ centered at a point $y\in  K$, there is $0\le j\le
N_0$ such that $\dist(f^{-j}(q),y)<\rho$.  Then, by the choice of
$\alpha $ and $\rho$, we have that  $W^s_{\delta_s}(f^{-j}(q))$
intersects $\gamma_1^u$ in a point $z$ with
$\dist_{W^s_{\delta_s}(f^{-j}(q))}(z,f^{-j}(q))<\alpha  <\delta_s/4$
and $\dist_{\gamma_1^u}(z,y)<\delta_1/2$. 
%
%
\cqd

\cle\label{l.delta2} There is $\delta_2>0$ such that if $\gamma^u$
is a centre-unstable disk of radius $\delta_1/2$ centered at a
point $z\in W^s_{\delta_s}(w)$ with $w\in  K$, then
$f^j(\gamma^u)$ contains a centre-unstable disk of radius
$\delta_2$ centered at $f^j(z)$, for each $1\le j\le N_0$. \fle

\dem 

Let us first prove the result for $j=1$. Let $f(y)$ be a point in
$\partial f(\gamma^u)$ minimizing the distance from $f(z)$ to
$\partial f(\gamma^u)$, and let $\eta_1$ be a curve of minimal
length in $f(\gamma^u)$ connecting $f(z)$ to $f(y)$. Define
$\eta_0=f^{-1}(\eta_1)$. 
Denote by $\dot\eta_1(x)$ the tangent vector to the curve $\eta_1$
at the point $x$. Then, 
$$
\|Df^{-1}(w) \dot\eta_1(x)\| \le C \, \|\dot\eta_1(x)\|,
$$
where
 $$C=\max_{x\in M}\big\{\|Df^{-1}(x)\|\big\}\ge1.$$
  Hence,
$$
\length(\eta_0) \le C\length(\eta_1) .
$$
Noting that $\eta_0$ is a curve connecting $z$ to $y\in \partial
\gamma^u$, this implies that $\length(\eta_0)\ge \delta_1/2$. Hence
$$
\length(\eta_1) \ge C^{-1}\length(\eta_0)\ge C^{-1}\delta_1/2.
$$
Thus  $f(\gamma^u)$ contains the  disk $\gamma^u_1$ of radius
$C^{-1}\delta_1/2$ around $f(z)$. Moreover,
 $$\dist(f(z),f(w))\le \lambda\delta_s<\delta_s,$$ and so, by the choice of
 $\delta_s$, we have that $\gamma^u_1$ is also a centre-unstable disk.
 Making now $\gamma^u_1$ play the role of $\gamma^u$ and $f^2(z)$ play the
 role of $f(z)$ we prove that:
 \begin{enumerate}\item[(a)] $f(\gamma^u_1)$
 contains a centre-unstable disk of radius $C^{-2}\delta_1/2^2$ centered at $f^2(z)$;
 \item[(b)] $\dist(f^2(z),f^2(w))\le
 \lambda^2\delta_s<\delta_s$.
 \end{enumerate}
 Item (a) gives in particular that  $f^2(\gamma^u)$
 contains a centre-unstable disk of radius $C^{-2}\delta_1/2^2$
 centered at
 $f^2(z)$.
 Arguing inductively we are able to prove that
 $f^j(\gamma^u)$ contains a disk of radius $C^{-j}\delta_1/2^j\ge C^{-N_0}\delta_1/2^{N_0}$
 around $f^j(z)$, for each $1\le j\le N_0$. Hence, we just have to take
 $\delta_2=C^{-N_0}\delta_1/2^{N_0}$.
 \cqd

\subsection{Partition on the reference leaf}\label{s.partion}

The construction that we are going to explain below requires the
use several constants. 
First we
take $\delta_1>0$ as in \eqref{e.delta1}, and
$0<\delta_2<\delta_1$ as in Lemma~\ref{l.delta2}. Then we take
$\delta_0>0$ and $\vare>0$ so that
 $$
\delta_0\ll\delta_2 \qand \varepsilon \ll \delta_0
 .
 $$
Next we describe the construction of the  ($ m_D$ mod 0)
partition $\cp$ of the unstable disk of radius $\delta_0$
centered at $p$ contained in $D$. For that we consider the following
 neighborhoods of \( p \) in $D$
$$\Delta^{0}_{0}=
B^u(p,\delta_{0}),\quad \Delta^{1}_{0}=B^u(p,2\delta_{0}),\quad
\Delta^{2}_{0}=B^u(p,\sqrt\delta_{0})\qand
\Delta^{3}_{0}=B^u(p,2\sqrt\delta_{0}),$$ and the cylinders over these sets,
 $$
 \cc^i=\bigcup_{x\in\Delta_0^i}W^s_{\delta_s}(x), \quad
 \text{for $i=0,1,2,3$.}$$
Letting $\pi$ denote the projection from $\cc^3$ onto $\Delta^{3}_{0}$
along local stable leaves, we have
 $$
 \pi(\cc^i)=\Delta^{i}_{0}, \quad
 \text{for $i=0,1,2,3$.}
 $$
We say that a centre-unstable disk $\gamma^u$ {\em $u$-crosses}
$\cc^i$ if $\pi(\gamma^u)=\Delta^{i}_{0}$.

\cre To simplify the exposition we shall pretend  that each
centre-unstable disk $\gamma^u$ $u$-crossing $\cc^i$ is still a
disk centered at a point in $W^s_{\delta_s}(p)$ with the same
radius of~$\Delta_0^i$. Actually, the radius of such a disk $\gamma^u$
is proportional to the
radius of $\Delta_0^i$, with the proportionality depending only on the height
of the cylinder and the angles of the two fibre bundles in the
dominated splitting. \fre

Let $\partial^u\cc_1^3$ denote the {\em top} and {\em bottom}
components of $\partial\cc^3$, i.e. the set of points
$z\in\partial\cc^3$ such that $z\in
\partial W^s_{\delta_s}(x)$ for some $x\in \Delta^{3}_{0}$. By the
domination property, we may take $\delta_0>0$ small so that no
centre-unstable disk contained in $\cc^3$ and intersecting
$W^s_{\delta_s/2}(p)$ can reach~$\partial^u\cc^3$.
For $0<\sigma<1$ given by Lemma~\ref{l:hyperbolic2}, let
\[
I_{k}=\left\{x\in\Delta^{1}_{0}\: : \:\delta_{0}(1+\sigma^{k/2}) <
\dist_D(x,p) < \delta_{0}(1+\sigma^{(k-1)/2})\right\},\quad  k\ge 1,
\]
be a partition ($\leb_D$ mod 0) into countably many rings of \(
\Delta_0^{1}\setminus \Delta_0^0 \).

 We are now able to start with the
 construction of the partition $\cp$ of $\Delta_0$. The
construction requires that we introduce inductively several objects.
In particular, we will consider sequences of sets $(\Delta_n)$,
$(A_n)$,  and $(B_n)$. For each  $n\ge0$, the set $\Delta_n$ is that
part of $\Delta_0$ that has not yet been partitioned up to time $n$. The set $\Delta_n$
is the disjoint union of $A_n$ and $B_n$, where $A_n$ is essentially
the part of $\Delta_n$ where new elements of partition may be
appear in the next step of the construction, and $B_n$ is some protection that we
put around the sets previously constructed in order to avoid
overlaps. For technical reasons, a small neighborhood $A_n^\vare$ of
each $A_n$ will also be considered.

\smallskip

\subsubsection*{First step of induction}
We fix \( R_{0} \) 
some large integer, and we ignore any dynamics occurring up to time
\( R_{0} \). Let \( k\geq R_{0}+1 \) be the first time that \(
\Delta_0\cap H_{k}\neq\emptyset \). For \(j < k
\) we define formally the objects 
 $$
A_{j}=A_{j}^{\varepsilon}=\Delta_{j}=\Delta_{0},\qand B_j=\emptyset
.
 $$
  Let
$(\omega_{k,j}^3)_j$ be all the center-unstable disks in
 $
 A_{k-1}^\vare
 $
contained in hyperbolic pre-balls $V_{m}$, with $k- N_0\le m\le k$,
which are mapped by $f^k$ onto a centre-unstable disk $u$-crossing
 $\cc^3$ and intersecting $W^s_{\delta_s/4}(p)$.
 Then we let
\[
\omega_{k,j}^i=\omega_{k,j}^3\cap f^{-k}(\cc^i),\quad i=0,1,2
\]
and set $R(x)=k$ for $x\in \omega_{k,j}^0$. We  take
\[
\Delta_{k}=\Delta_{k-1}\setminus \{R=k\},
\]
and define a function \( t_{k}:\Delta_{k}\to \mathbb N \) by
\begin{equation*}
    t_{k}(x) =
    \begin{cases}
    s, & \text{ if }  x\in \omega_{k,j}^{1} \text{ and }
    \pi(f^{k}(x)) \in I_{s} \text{  for some $j$;} \\
    0 ,& \text{ otherwise}.
    \end{cases}
\end{equation*}
Finally let
\[
A_{k}= \{x\in\Delta_{k}: t_{k}(x) = 0\}, \quad B_{k}=
\{x\in\Delta_{k}: t_{k}(x) > 0\}
\]
and
 $$
 A_{k}^\vare= \{x\in\Delta_{k}:
 \dist_{f^{k+1}(D)}(f^{k+1}(x),f^{k+1}(A_k))<\vare\}.
 $$

\subsubsection*{General step of induction}

The general step of the construction follows the ideas above with
minor modifications.  Assume that the sets $\Delta_{i}$, $A_{i}$,
$A_{i}^\vare$
 $B_{i}$, $\{R=i\}$ and functions \( t_{i}:
\Delta_{i}\to\mathbb N \) are defined for each \( i\leq n-1 \). 
 Let
$(\omega_{n,j}^3)_j$ be all the centre-unstable disks in
 $
 A_{n-1}^\vare
 $
contained in hyperbolic pre-balls $V_{m}$, with $n- N_0\le m\le n$,
which are mapped by~$f^n$ onto a centre-unstable disk $u$-crossing
 $\cc^3$ and intersecting $W^s_{\delta_s/4}(p)$.
 Take
\begin{equation}\label{eqomega}
\omega_{n,j}^i=\omega_{n,j}^3\cap f^{-n}(\cc_1^i),\quad i=0,1,2
\end{equation}
set $R(x)=n$ for $x\in \omega_{n,j}^0$, and let
\[
\Delta_{n}=\Delta_{n-1}\setminus \{R=n\}.
\]
The definition of the function \( t_{n}:\Delta_{n}\to \mathbb N \)
is slightly different in the general case:
\begin{equation*}
    t_{n}(x) =
    \begin{cases}
    x ,& \text{ if } x\in \omega_{n,j}^{1}\setminus \omega_{n,j}^{0} \text{ and }
    f^{n}(x) ,\in I_{s} \text{ for some $j$,} \\
    0, & \text{ if } x\in A_{n-1} \setminus \bigcup_{j} \omega^{1}_{n,j},\\
    t_{n-1}(x)-1 ,& \text{ if } x\in B_{n-1}\setminus \bigcup_{j} \omega^{1}_{n,j}.
    \end{cases}
\end{equation*}
Finally, we let
\[
A_{n}= \{x\in\Delta_{n}: t_{n}(x) = 0\}, \quad B_{n}=
\{x\in\Delta_{n}: t_{n}(x) > 0\}
\]
and
 $$
 A_{n}^\vare= \{x\in\Delta_{n}:
 \dist_{f^{n+1}(D)}(f^{n+1}(x),f^{n+1}(A_n))<\vare\}.
 $$
At this point we have described the construction of the sets  $A_n$,
$A_n^\vare$,  $B_n$ and $\{R=n\}$.

\medskip

Since the components of  $\{R=n\}$ are taken in $A_{n-1}^\vare$,
it could happen that these new components intersect $B_{n-1}$. The
next lemma shows that this is not the case as long as $\epsilon>0$ is taken small enough.
For notational simplicity we will drop the index $j$ in the
elements defined at \eqref{eqomega}.

\cle\label{l.claim} If $\epsilon>0$ is small, then
$\omega_{n}^1\cap\{t_{n-1}\ge 1\}=\emptyset$ for all $n\ge 1$. \fle

\dem Take  $k\ge1$ and let  $\omega_{n-k}^0$ be a component of
$\{R=n-k\}$. Let $Q_{k}$ be the part of $\omega_{n-k}^1$ that is
mapped by $\pi\circ f^{n-k}$ onto $I_{k}$, and assume that $Q_{k}$
intersects some $\omega_{n}^3$. Recall that, by construction,
$Q_{k}$ is precisely that part of $\omega_{n-k}^1$ on which
$t_{n-1}=1$, and $\omega_{n}^3$ is contained in a hyperbolic
pre-ball $V_{m}$ with $n-N_0\le m\le n$.

Let $q_1$ and $q_2$ be any two points in distinct components (inner
and outer, respectively) of the boundary of $Q_{k}$. If we assume
that $q_1,q_2\in \omega_{n}^3$, then $q_1,q_2\in V_m$, and so by
Lemma~\ref{l.contraction} we have
 \begin{equation}\label{e.zq1}
 \dist_{f^{n-k}(D)}(f^{n-k}(q_1),f^{n-k}(q_2))\le
 C_0\sigma^{k/2}\dist_{f^n(D)}(f^{n}(q_1),f^{n}(q_2))
 \end{equation}
for some $C_0$ depending on $N_0$. We also have for some $C_1>0$
depending on the angle of the stable and centre-unstable spaces over
$K_\infty$
 \begin{eqnarray*}
\dist_{f^{n-k}(D)}(f^{n-k}(q_1),f^{n-k}(q_2))&\ge& C_1
\delta_{0}(1+\sigma^{(k-1)/2})-\delta_{0}(1+\sigma^{k/2})\\&=& C_1
\delta_{0} \sigma^{k/2}(\sigma^{-1/2}-1),
 \end{eqnarray*}
which combined with (\ref{e.zq1})  gives
$$
 \dist_{f^n(D)}(f^{n}(q_1),f^{n}(q_2))\ge \frac{C_1}{C_0}\delta_{0}(\sigma^{-1/2}-1).
 $$
On the other hand, since $ \omega^{3}_{n}\subset
A_{n-1}^{\varepsilon}$ by construction of $\omega^{3}_{n}$, taking
 \begin{equation}\label{eq.vare}
 \vare<\frac{C_1}{C_0}\delta_{0}(\sigma^{-1/2}-1)
 \end{equation}
we  have $\omega_{n}^3\cap\{t_{n-1}>1\}=\emptyset$. This implies
$\omega_{n}^1\cap\{t_{n-1}\ge 1\}=\emptyset$.\cqd

 \subsection{Product structure}


Consider the center-unstable disk $\Delta_0\subset D$ and the
partition $\cp$ of $\Delta_0$ ($\leb_D$ mod 0) defined in
Section~\ref{s.partion}. We shall use the elements of $\cp$ to
define the $s$-subsets that give rise to the hyperbolic structure.
Given an arbitrary element $\omega\in\cp$, we have by construction
some $R(\omega)\in\NN$ such that $f^{R(\omega)}(\omega)$ is a
centre-unstable disk $u$-crossing~$\cc^0$. We define $\cc_\omega$ as
the cylinder made by the stable leaves passing through the points in
$\omega$, i.e.
 $$\cc_\omega=\bigcup_{x\in\omega}W^s_{\delta_s}(x).$$
The sets $\cc_\omega$, with $\omega\in\cp$, are by definition the
pairwise disjoint $s$-subsets $\Lambda_1,\Lambda_2,\dots$ which
define the Markovian structure.

Now we define inductively some sets of centre-unstable manifolds
$u$-crossing $ \cc^0$ that will give rise to the family  $\Gamma^u$.
The first one is
 $$
 \Gamma_0=\left\{\Delta_0\right\}.
 $$
Having defined  $\Gamma_j$,  for some $j\ge0$, we define
 $$
 \Gamma_{j+1}=\left\{f^{R(\omega)}(\cc_\omega\cap \gamma)\,\colon\,
 \,
 \omega\in\cp \text{ and } \gamma\in\Gamma_j\right\}.
 $$
 Observe that each element of $\Gamma_j$ is equal to 
 an iterate of a subset of $\Delta_0$. In
 particular, the elements of each $\Gamma_j$ are unstable
 manifolds. Moreover, since by construction $f^{R(\omega)}(\omega)$ intersects
 $W^s_{\delta_s/4}(p)$, then according to the choice of
 $\delta_0$ and the invariance of the stable foliation, we have
 that each element of $\Gamma_j$ must $u$-cross $\cc^0$.


Since the union of the leaves of the sets $\Gamma_j$, with
$j\ge0$, is not necessarily compact, we still need to take
accumulation points of that union. Let
 $$
\Delta_\infty=\overline{\bigcup_{j\ge0}\bigcup_{\gamma_j\in\Gamma_j}\gamma_j}.
 $$
Given $x\in\Delta_\infty$, there are $(j_k)_k\to\infty$,  disks
$\gamma_{j_k}\in \Gamma_{j_k}$ and points $x_k\in\gamma_{j_k}$
converging to $x$ as $k\to\infty$. Using the domination property and
Ascoli-Arzela theorem we conclude that the disks $\gamma_{j_k}$
converge to a disk $\gamma_\infty$  containing $x$. Since the disk
$\gamma_\infty$ is accumulated by disks $u$-crossing $\cc^0$ then it
also must $u$-cross $\cc^0$. We define $\Gamma_\infty$ as the set of
all these accumulation disks.
 Finally, we take
  $$\Gamma^u=\bigcup_{j\ge0}\Gamma_j\cup\Gamma_\infty.$$

  \subsection{Backward contraction and bounded distortion}
The backward contraction property (P$\bf_3$) follows from Lemma~\ref{le.contravolta} below. Bounded distortion is typically a consequence of backward contraction together with some Hölder control of $\log|\det Df^u|$. Property (P$\bf_4$) follows naturally from Proposition~\ref{c.curvature}  together with Lemma~\ref{le.contravolta}  exactly as in the proof of \cite[Proposition~2.8]{ABV}.

\cle\label{le.contravolta} There is $C>0$ such that, given   $\omega\in\cp$  and
$\gamma\in\Gamma^u$, we have for all $1\le k \le R(\omega)$ and all
$x,y\in \cc_\omega\cap\gamma$
 $$
\dist_{f^{R(\omega)-k}(\cc_\omega\cap\gamma)}(f^{R(\omega)-k}(x),f^{R(\omega)-k}(y))
\le
C\sigma^{k/2}\dist_{f^{R(\omega)}(\cc_\omega\cap\gamma)}(f^{R(\omega)}(x),f^{R(\omega)}(y)).
 $$
 \fle

\dem 

Recall that, by construction, for each $\omega\in\cp$  there is a
hyperbolic pre-ball $V_{n(\omega)}(x)$ containing $\omega$
associated to some point $x\in D$ with $\sigma$-hyperbolic time
$n(\omega)$ satisfying $R(\omega)-N_0\le n(\omega)\le R(\omega)$.
Taking $\delta_s,\delta_0<\delta_1/2$,  it follows from
\eqref{e.delta1} that $n(\omega)$ is  a $\sqrt\sigma$-hyperbolic
time for every point in $\cc_\omega\cap\gamma$. Then, recall
\eqref{contra}, this implies that for all $1\le k \le n(\omega)$
and all $x,y\in \cc_\omega\cap\gamma$ we have
 $$
\dist_{f^{n(\omega)-k}(\cc_\omega\cap\gamma)}(f^{n(\omega)-k}(x),f^{n(\omega)-k}(y))
\le
\sigma^{k/2}\dist_{f^{n(\omega)}(\cc_\omega\cap\gamma)}(f^{n(\omega)}(x),f^{n(\omega)}(y)).
 $$
Since the difference between $R(\omega)$ and $n(\omega)$ is at most
$N_0$, the  result follows with $C$ depending only on $N_0$ and the
the derivative of $f$. \cqd

\subsection{Regularity of the foliations}\label{sec.regularity}


Here we prove property (P$_5$). 
This is standard for uniformly hyperbolic attractors, and we shall adapt classical ideas to our setting. We begin with the statement of a useful result on vector bundles whose proof can be found in \cite[Theorem 6.1]{HP}.


\cle\label{th.hirsch-pugh}
Let $p\colon Y\to X$ be a vector bundle over a metric space $X$ endowed with an admissible metric. Let $D\subset Y$ be the unit ball bundle, and $F\colon D\to D$ a map covering a homeomorphism $f\colon X\to X$. Suppose $0\le \kappa<1$ and that for each $x\in X$, the restriction $F_x\colon D_x\to D_x$ has Lipschitz constant $\lip(F_x)\le \kappa$. Then
\begin{enumerate}
  \item There is a unique section $\sigma_0\colon X\to D $ whose image is invariant under $F$.
  \item Let $\lip(f)=c<\infty$ and $0<\alpha\le 1$ be such that $\kappa c^\alpha<1$. Then $\sigma_0$ satisfies a Hölder condition of exponent $\alpha$.
\end{enumerate}
\fle

For the sake of completeness, let us mention that a metric $d$ on $E$ is \emph{admissible} if there is a complementary bundle $E'$ over $X$, and an isomorphism $h\colon E\oplus E'\to X\times A$ to a product bundle, where $A$ is a Banach space, such that $d$ is induced from the product metric on $X\times A$.

\ct\label{th.Holder}
Let $f: M\to M$ be a $C^1$ diffeomorphism and let $K\subset M$ be a compact  invariant  set with a  dominated splitting   $T_K M=E^{cs}\oplus
E^{cu}$. Then the fiber bundles $E^{cs}$ and $E^{cu}$ are Hölder continuous on $K$.
\ft

\dem We consider the centre-unstable bundle, the other one is similar.
For each $x\in K$ let $L_x$ be the space of bounded linear maps $L(E_x^{cu},E_x^{cs})$. For each $x\in K$, let $L_x(1)$ denote the unit ball around $0\in L_x$, and define  $\Gamma_x: L_x(1)\to L_{f(x)}(1)$ as the graph transform induced by
 $Df(x)$: 
    $$\Gamma_x(\mu_x)=(Df \vert E^{cs}_x)\cdot\mu_x \cdot (Df^{-1} \vert E^{cu}_{f(x)}).$$
Consider $L(E^{cu}, E^{cs})$ the vector bundle over $K$ whose fiber over $x\in K$ is $L_x$, and let $D$ be its unit ball bundle. Then $\Gamma : D\to D$ is a bundle map covering $f\vert K$ with
 $$\lip(\Gamma_x)\le \|Df \mid E^{cs}_x\|
\cdot \|Df^{-1} \mid E^{cu}_{f(x)}\| \le\lambda<1.$$
Let $c$ be a Lipschitz constant for $\vert K$, and choose $0<\alpha\le 1$ small so that  $\lambda c^\alpha<1$. By Lemma~\ref{th.hirsch-pugh} there exists a unique section $\sigma_0\colon X\to D $ whose image is invariant under~$F$ and it satisfies a Hölder condition of exponent~$\alpha$. This unique section is necessarily the null section.
\cqd

The next result gives precisely (P$\bf_5$)(a).

\cco\label{co.produtorio}
There are $C>0$ and $0<\beta<1$ such that for all $y\in\gamma^s(x)$ and $ n\ge 0$
 $$\displaystyle
 \log \prod_{i=n}^\infty\frac{\det Df^u(f^i(x))}{\det Df^u(f^i(y))}\le C\beta^{n}.
 $$
\fco
\dem As we are assuming that $Df$ is Hölder continuous, it follows from Theorem~\ref{th.Holder} that $\log|\det Df^u|$ is  Hölder continuous. The conclusion is then an immediate consequence of the uniform contraction on stable leaves.
\cqd

Now we are going to prove (P$\bf_5$)(b). We start by introducing some useful notions. We say that $\phi: N\to P$, where $N$ and $P$ are submanifolds of $M$,  is \emph{absolutely continuous} if it is an injective map for which there exists $J:N\to\RR$, called the \emph{Jacobian} of $\phi$, such that
 $$
 \leb_P(\phi(A))=\int_A Jd\leb_N.
 $$
Property (P$\bf_5$)(b) can be restated in the following terms:

\cpr\label{pr.regulstable}
Given
$\gamma,\gamma'\in\Gamma^u$, define
$\phi\colon\gamma'\to\gamma$  by
$\phi(x)=\gamma^s(x)\cap \gamma$. Then $\phi$ is absolutely continuous and  the  Jacobian  of $\phi$ is given by
        $$
        J(x)=
        \prod_{i=0}^\infty\frac{\det Df^u(f^i(x))}{\det
        Df^u(f^i(\phi(x)))}.$$
\fpr
One can easily deduce from Corollary~\ref{co.produtorio} that this infinite product converges uniformly.
The remaining of this section is devoted to the proof of Proposition~\ref{pr.regulstable}.
We start with a general result about the convergence of Jacobians whose proof is given in \cite[Theorem~3.3]{Ma}.

\cle\label{le.jacomane}
Let $N$ and $P$ be manifolds, $P$ with finite volume, and for each $n\ge 1$,  $\phi_n:N\to P$ an
absolutely continuous map with Jacobian $J_n$. Assume that
\begin{enumerate}
  \item $\phi_n$ converges uniformly to an injective
continuous map $\phi:N\to P$;
  \item $J_n$ converges uniformly to an integrable function $J: N\to\RR$.
\end{enumerate}
Then
$\phi$ is absolutely continuous with Jacobian $J$.
\fle
For the sake of completeness, let us mention that there is a slight difference in our definition of absolute continuity of maps. Contrarily to  \cite{Ma}, and for reasons that will become clear in later, we do not impose continuity of the maps $\phi_n$. However, the proof of \cite[Theorem~3.3]{Ma} uses only the continuity of the limit function $\phi$, and so it still works in our case.

Consider now $\gamma,\gamma'\in\Gamma^u$ and
$\phi\colon\gamma'\to\gamma$  as in Proposition~\ref{pr.regulstable}. The proof of the  next lemma is given in \cite[Lemma 3.4]{Ma} for uniformly hyperbolic diffeomorphisms. Nevertheless, one can easily see that it  is obtained as a consequence of \cite[Lemma 3.8]{Ma} whose proof uses only the existence of a dominated splitting.

\cle\label{le.mane}
For each $n\ge 1$, there is an absolutely continuous  $\pi_n:f^n(\gamma)\to f^n(\gamma')$ with Jacobian $G_n$ satisfying
\begin{enumerate}
  \item $\displaystyle \lim_{n\to\infty}\sup_{x\in \gamma}\left\{ \dist_{f^n(\gamma')}(\pi_n(f^n(x)),f^n(\phi(x))\right\}=0$;
  \item $\displaystyle \lim_{n\to\infty}\sup_{x\in f^n(\gamma)} \left\{ |1-G_n(x)|  \right\}=0$.
\end{enumerate}

\fle

We consider the sequence of  consecutive return times
for points in $\Lambda$,
 \begin{equation*}\label{def.rs}
 s_1=R\qand s_{n+1}=s_{n}+R\circ f^{s_{n}},\quad\text{for $n\ge1$}.
 \end{equation*}
Notice that these return time functions are defined $\leb_\gamma$ almost everywhere on each $\gamma\in \Gamma^u$ and are piecewise constant.

\cre\label{re.sp} Using the sequence of return times one can easily construct a sequence of  ($\leb_\gamma$ mod 0 on each $\gamma\in \Gamma^u$) partitions   $(\cp_n)_n$ by $s$-subsets of $\Lambda$ with $s_n$ constant on each element of $\cp_n$, for which (P$_1$)-(P$_5$) hold when we take $s_n$ playing the role of $R$ and the elements of $\cp_n$ playing the role of the $s$-subsets. Moreover, the constants $C>0$ and $0<\beta<1$ can be chosen not depending on $n$.
\fre

We define, for each $n\ge 1$, the map $\phi_n: \gamma\to\gamma'$ as
 \begin{equation}\label{eq.phin}
 \phi_n=f^{-s_n}\pi_{s_n} f^{s_n}.
 \end{equation}
It is straightforward to check that $\phi_n$ is absolutely continuous with   Jacobian
  \begin{equation}\label{eq.Jn}
  J_n(x)=\frac{|\det (Df^{s_n})^u(x)|}{|\det
        (Df^{s_n})^u(\phi_n(x))|}\cdot G_{s_n}(f^{s_n}(x)).
  \end{equation}
Observe that these functions are defined $\leb_\gamma$ almost everywhere. So, we may find a Borel set $F\subset \gamma$ with full $\leb_\gamma$ measure on which they are all defined. We extend $\phi_n$ to $\gamma$ simply by considering $\phi_n(x)=\phi(x)$ and $J_n(x)=J(x)$ for all $n\ge 1$ and $x\in \gamma\setminus F$. Since $F$ has zero $\leb_\gamma$ measure one still has that $J_n$ is the Jacobian of $\phi_n$.

Proposition~\ref{pr.regulstable} is and immediate consequence of  Lemma~\ref{le.jacomane} together with the next~one.

\cle\label{le.phin}
The maps $\phi_n$ converge   uniformly to $\phi$, and  the  Jacobians $J_n$ converge uniformly to $J$.
\fle

\dem It is enough to prove the convergence of each sequence restricted to $F$, i.e restricted to the set of points where the expressions of $\phi_n$ and $J_n$ are given by~\eqref{eq.phin} and~\eqref{eq.Jn} respectively.

Let us prove first  that $\phi_n$ converges uniformly to $\phi$. Using the backward contraction on unstable leaves given by (P$\bf_3$), recall Remark~\ref{re.sp}, we may write for $x\in\gamma$
%
\begin{eqnarray*}
  \dist_{\gamma'}(\phi_n(x),\phi(x))  &=& \dist_{\gamma'}(f^{-s_n}\pi_{s_n} f^{s_n}(x),f^{-{s_n}}f^{s_n}\phi(x)) \\
   &\le& C\beta^{s_n}\dist_{f^{s_n}(\gamma')}(\pi_{s_n} f^{s_n}(x),f^{s_n}\phi(x)).
\end{eqnarray*}
Since ${s_n}\to\infty$ as $n\to\infty$ and $\dist_{f^{s_n}(\gamma')}(\pi_{s_n} f^{s_n}(x),f^{s_n}\phi(x))$ is bounded by Lemma~\ref{le.mane}, we have the uniform convergence of $\phi_n$ to $\phi$.

 By \eqref{eq.Jn} we have
  $$J_n(x)=\frac{|\det (Df^{{s_n}})^u(x)|}{|\det
        (Df^{{s_n}})^u(\phi (x))|}\cdot \frac{|\det (Df^{{s_n}})^u(\phi (x))|}{|\det
        (Df^{{s_n}})^u(\phi_n(x))|}\cdot G_{{s_n}}(f^{{s_n}}(x)).
        $$
 Using the chain rule and  Corollary~\ref{co.produtorio} it easily follows that the first term in the product above converges uniformly to $J(x)$. Moreover, the third term converges uniformly to one, by Lemma~\ref{le.mane}. It remains to see that the middle term also converges uniformly to one. Recalling Remark~\ref{re.sp}, by  bounded distortion we have
  \begin{eqnarray*}
    \frac{|\det (Df^{{s_n}})^u(\phi (x))|}{|\det
        (Df^{{s_n}})^u(\phi_n(x))|}  &\le &\exp\big(C\dist_{f^{s_n}(\gamma')}(f^{s_n}(\phi(x)),f^{s_n}(\phi_n(x)))^{\eta}\big) \\
     &=& \exp\big(C\dist_{f^{s_n}(\gamma')}(f^{s_n}(\phi(x)),\pi_{s_n}(f^{s_n}(x)))^{\eta}\big).
  \end{eqnarray*}
  Similarly we obtain
  $$\frac{|\det (Df^{{s_n}})^u(\phi (x))|}{|\det
        (Df^{{s_n}})^u(\phi_n(x))|}  \ge  \exp\big(-C\dist_{f^{s_n}(\gamma')}(f^{s_n}(\phi(x)),\pi_{s_n}(f^{s_n}(x)))^{\eta}\big).
        $$
        The conclusion then follows from Lemma~\ref{le.mane}.
 \cqd

\section{Decay estimates}\label{s.metric}

In this section we obtain the metric estimates on the decay of $\leb_D\{R> n\}$ of
Theorem~\ref{t:Markov towers}. These estimates are an adaptation of similar ones from \cite[Section~5]{ALP} to our setting. We start by proving estimates
arising directly from the construction preformed in
Section~\ref{s.partion}. In the final part of the argument we use some results that have been put into an abstract setting in  \cite[Section~4.5.2]{Al} and get the desired
conclusion.

 \cle\label{l.flowb} There exists $a_0>0$
  such that  $$\leb_D(B_{n-1}\cap A_n)\ge a_0 \leb_D(B_{n-1}),$$ for
  all
$n\ge1$.
 \fle
  \dem
It is enough to see this for each component of $B_{n-1}$. Let  $C$
be a component of $B_{n-1}$ and let $Q$ be  its outer ring,
corresponding to $t_{n-1}=1$. Observe that by Lemma~\ref{l.claim} we
have $Q=C\cap A_n$. Moreover, there must be some $k<n$ and a
component $\omega^0_k$ of $\{R=k\}$ such that $\pi\circ f^k$ maps
$C$ diffeomorphically onto $\cup_{i=k}^\infty I_i$ and $Q$ onto
$I_k$, both with uniform bounded distortion as in Corollary~\ref{p.distortion}. 
 Thus, it is sufficient to compare the Lebesgue measures of $\cup_{i=k}^\infty
I_i$ and $ I_k$. We have
 $$
 \frac{ \leb_D( I_k)}{\leb_D(\cup_{i=k}^\infty
I_i)}\thickapprox\frac{
[\delta_0(1+\sigma^{(k-1)/2})]^u-[\delta_0(1+\sigma^{k/2})]^u}
{[\delta_0(1+\sigma^{(k-1)/2})]^u-\delta_0^u}\thickapprox
1-\sigma^{1/2},
 $$
where $u$ is the dimension of $E^{cu}$.
  \cqd

  \cle\label{l.flowa}  There exist $b_0,c_0>0$  with 
  $b_0,c_0\to 0$  as $\delta_0\to 0$,
  such that
\begin{enumerate}
 \item $\leb_D(A_{n-1}\cap B_n)\le b_0 \leb_D(A_{n-1})$, \item
$ \leb_D(A_{n-1}\cap \{R=n\})\le c_0 \leb_D(A_{n-1})$,
\end{enumerate}
for all $n\ge1$.
\fle \dem It is enough to prove this for each  neighborhood of a
component $\omega^0_n$ of $\{R=n\}$. Observe that by construction we
have $\omega^3_n\subset A_{n-1}^\vare$, which means that $\omega^2_n
\subset A_{n-1}$, because we are taking $\vare<\delta_0$. Using the
uniform bounded distortion of Corollary~\ref{p.distortion} we obtain
 $$
 \frac{ \leb_D(\omega^1_n\setminus \omega^0_n)}{\leb_D(\omega^2_n\setminus \omega^1_n)}
 \thickapprox
 \frac{ \leb_D(\Delta^1_0\setminus \Delta^0_0)}{ \leb_D(\Delta^2_0\setminus \Delta^1_0)}
 \thickapprox
 \frac{\delta_0^d}{\delta_0^{d/2}}\ll 1,
 $$
which gives the first estimate. Moreover,
$$
 \frac{ \leb_D( \omega^0_n)}{\leb_D(\omega^2_n\setminus \omega^1_n)}
 \thickapprox
 \frac{ \leb_D( \Delta^0_0)}{\leb_D(\Delta^2_0\setminus \Delta^1_0)}
 \thickapprox
 \frac{\delta_0^d}{\delta_0^{d/2}}\ll 1,
$$
and this gives the second one.
 \cqd

    \cpr\label{p.construction}
    There exist $c_1>0$ and a positive integer $N=N(\vare)$ such that
    \[
      \leb_D\left(\cup_{i=0}^N\big\{R=n+i\big\}\right)\ge c_1  \leb_D(A_{n-1}\cap H_{n}),
    \] for all \( n\ge1\).
    \fpr
\dem Let $K_0=\max_{x\in\Lambda}\|Df^{-1}\|$ and take
$r=5\delta_0K_0^{N_0}$, where $N_0$ is given by Lemma~\ref{l.N0q}.
 Recall that by Lemma~\ref{l.contraction}, for
each $z\in f^n(A_{n-1}\cap H_n)$ there is $x\in H_n$ and a
$\sigma$-hyperbolic pre-ball $V_n(x)\subset D$ which is sent
diffeomorphically onto the centre-unstable ball of radius
$\delta_1$ around $z$. Let $\{z_j\}$ be a maximal set in
$f^n(A_{n-1}\cap H_n)$ with the property that the sets $B_r(z_j)$
are pairwise disjoint, where each $B_r(z_j)$  is the ball of
radius $r$ centered at $z_j$ inside the hyperbolic ball around
$z_j$. By maximality we have
 \begin{equation}\label{eq.cont}
 \mcup_j B_{2r}(z_j)\supset f^n(A_{n-1} \cap H_n).
 \end{equation}
For each $j$ let $x_j\in H_n$  be the point such that
$f^n(x_j)=z_j$.

\smallskip

\noindent {\sc Claim 1.} {\em There is  $0\le k\le N_\epsilon+N_0$
such that $t_{n+k}$ is not identically zero in
$f^{-n}(B_\epsilon(z))$.}

\smallskip \noindent
  Assume, by contradiction, that $ t_{n+k}\vert f^{-n}(B_\epsilon(z))=0$ for all $0\le
k\le N_\epsilon+N_0 $. This implies that $f^{-n}(B_\epsilon(z))
\subset A_{n+k}$ for all $0\le k\le N_\epsilon+N_0$. Using
Lemma~\ref{le:grow1} we may find  a hyperbolic pre-ball
$V_{m}\subset B_\epsilon(z) $ with $\sigma$-hyperbolic time $m\leq
N_{\epsilon}$. Now, since $f^m(V_{m})$ is a centre-unstable disk
of radius $\delta_1$, it follows from Lemma~\ref{l.N0q} and
Lemma~\ref{l.delta2} that there are $V\subset f^m(V_{m})$ and
$m'\le N_0$ such that $u$-crossing $\cc^3$ and intersecting
$W^s_{\delta_s/4}(p)$. Thus, taking $k=m+m'$ we have that $0\le
k\le N_\epsilon+N_0$ and $f^{-n}(V_m)$ contains an element of
$\{R=n+k\}$ inside $f^{-n}(B_\epsilon(z))$. This contradicts the
fact that $t_{n+k}\vert f^{-n}(B_\epsilon(z))=0$ for all $0\le
k\le N_\epsilon+N_0$.

\smallskip

\noindent {\sc Claim 2.} {\em $f^{-n}(B_{\delta_1/4}(z))$ contains
a component of $\{R=n+k\}$ with $0\le k\le N_\epsilon+N_0$.}

\smallskip \noindent Let $k$ be the smallest integer $0\le k\le
N_\epsilon+N_0$ for which $t_{n+k}$ is not identically zero in $
f^{-n}(B_\epsilon(z))$. Since $f^{-n}(B_\epsilon(z))\subset
A_{n-1}^\epsilon\subset \{t_{n-1}\le1 \},$
 there must be some component
$\omega^{0}_{n+k}$ of $\{R=n+k\}$ for which
 $
 f^{-n}(B_\epsilon(z))\cap \omega_{n+k}^1\neq\emptyset.
 $
Recall that, by definition, $f^{n+k}$ sends $\omega_{n+k}^1$
diffeomorphically onto a centre-unstable disk (of radius
$2\delta_0$) $u$-crossing $\cc^1$ and intersecting
$W^s_{\delta_s/4}(p)$.  Thus, the diameter of
$f^n(\omega_{n+k}^1)$ is at most $4\delta_0K_0^{N_0}$. Since
$B_\epsilon(z)$ intersects $f^n(\omega_{n+k}^1)$ and
$\epsilon<\delta_0<\delta_0K_0^{N_0}$, we have
$f^{-n}(B_{\delta_1/4}(z))\supset \omega_{n+k}^0$,
 as long as we take $\delta_0>0$ small so that
 $
 5\delta_0K_0^{N_0}<\delta_1/4.
 $
Hence, we have shown  that $f^{-n}(B_{\delta_1/4}(z))$ contains
some component of $\{R=n+k\}$ with  $0\le k\le N_\epsilon+N_0$,
and so we have proved the claim.

\smallskip

Since $n$ is a hyperbolic time for $x_j$, we have by the
distortion control given by Corollary~\ref{p.distortion} that
there is some constant $C$ only depending on $C_2$ and $\delta_1$
for which
 \begin{equation}\label{eq.quo1}
 \frac{ \leb_D(f^{-n}(B_{2r}(z_j)))}{ \leb_D(f^{-n}(B_r(z_j)))}
 \le
 {C}\frac{ \leb_{f^n(D)}(B_{2r}(z_j))}{ \leb_{f^n(D)}(B_r(z_j))}
 \end{equation}
and
 \begin{equation}
\label{eq.quo2} \frac{ \leb_D(f^{-n}(B_r(z_j)))}{ \leb_D(
\omega_{n+k}^0)}
 \le
 {C}\frac{ \leb_{f^n(D)}(B_r(z_j))}{ \leb_{f^n(D)}(f^n(\omega_{n+k}^0))}.
 \end{equation}
Recalling that from time $n$ up to $n+k$ we have at most $N_0$
iterates, from \eqref{eq.quo1} and \eqref{eq.quo1} we easily
deduce that there is some positive constant, that we still denote
by $C$, for which
  $$
  \leb_D(f^{-n}(B_{2r}(z_j)))\le C \leb_D(
f^{-n}(B_r(z_j)))
 $$
 and $$
  \leb_D(f^{-n}(B_r(z_j)))\le C \leb_D(
\omega_{n+k_j}^0).
 $$
 Finally, let us compare the Lebesgue measure of
the sets $\bigcup_{i=0}^N\big\{R=n+i\big\}$ and $A_{n-1}\cap
H_{n}$. By \eqref{eq.cont} we have
 $$
 \leb_D\big(A_{n-1}\cap H_n \big)\le  \sum_j  \leb_D(f^{-n}(B_{2r}(z_j)))
 \le C
 \sum_j  \leb_D(f^{-n}(B_r(z_j))).
 $$
On the other hand, by the disjointness of the balls $B_r(z_j)$ we
have
$$
 \sum_j  \leb_D(f^{-n}(B_r(z_j)))\le C
 \sum_j  \leb_D( \omega_{n+k}^0) \le C
  \leb_D\left(\cup_{i=0}^N\big\{R=n+i\big\}\right).
  $$
We just have to take $c_1=C^{-2}$. \cqd

For completing the proof of Theorem~\ref{t:Markov towers}, it is
enough to show that
 $$
\leb_D\{\mathcal E> n\} ={O}(n^{-\tau})
\quad\Rightarrow\quad \leb_D\{R>n\}= {O}(n^{-\tau}).
 $$
Recall that  we have defined $H_n$, for $n\ge 1$, as the set of
points for which $n$ is a  $\sigma$-hyperbolic time. In
Corollary~\ref{c:hyperbolic3} we obtained the following
estimate:
  \begin{enumerate}
   \item[ (m$_1$) ] \emph{There is $\theta>0$ such that for all $n\ge1$ and
$A\subset M\setminus \{\mathcal E> n\}$ with $\leb_D(A)>0$}
    \[
   \frac{1}{n}\sum_{j=1}^{n}\frac{ \leb_D(A\cap H_{j})}{ \leb_D(A)} \geq \theta.
   \]
\end{enumerate}
   In the construction of the Markov structure have taken a disk $\Delta$ of radius $\delta_0>0$
   and defined inductively the subsets
  $A_n$, $B_n$, $\{R=n\}$ and $\Delta_n$ related in
   the following way:
    $$\Delta_n=\Delta\setminus\{R\le n\}=A_n\dot\cup B_n.$$
Moreover, we have proved in Lemma~\ref{l.flowb}, Lemma~\ref{l.flowa}
and Proposition~\ref{p.construction} that the following metric
relations hold:
  \begin{enumerate}
\item[(m$_2$) ] \emph{There is $a_0>0$ (bounded away from 0 for all
$\delta_0$) such that for all $n\ge1$}
 $${ \leb_D(B_{n-1}\cap A_n)}\ge a_0{ \leb_D(B_{n-1})}.$$
\item[(m$_3$) ]  \emph{There are $b_0,c_0>0$  with 
$b_0,c_0\to 0$
 as $\delta_0\to 0$, such
that for all $n\ge1$
 $$\frac{ \leb_D(A_{n-1}\cap B_n)}{ \leb_D(A_{n-1})}\le b_0\qand
 \frac{ \leb_D(A_{n-1}\cap \{R=n\})}{ \leb_D(A_{n-1})}\le c_0.$$}
\item[(m$_4$) ] \emph{There is $r_0>0$ and an integer $N\ge 0$ such that
for all \( n\ge1\)}
  $$  \leb_D\left(\cup_{i=0}^N\big\{R=n+i\big\}\right)\ge r_0 \leb_D(A_{n-1}\cap
    H_{n}).$$
\end{enumerate}
Estimates (m$_1$)-(m$_4$)  are enough to use the results of
\cite[Section 4.5.2]{Al} and obtain the  decay of
$\leb_\gamma\{R> n\}$ as in the conclusion of Theorem~\ref{t:Markov towers}.

\subsection*{Acknowledgment} The authors acknowledge the referee for valuable comments and references.


\begin{thebibliography}{10}

\bibitem{Aar} J. Aaronson, \emph{An introduction to infinite ergodic theory}, Mathematical Surveys and Monographs 50, American Mathematical Society, Providence, RI, 1997.

\bibitem{AD01} J. Aaronson, M. Denker, \emph{Local limit theorems for partial sums of stationary sequences generated by Gibbs-Markov maps},
Stoch. Dyn.  \textbf{1},  no. 2 (2001), 193-237.

\bibitem{AD02} J. Aaronson, M. Denker, \emph{Group extensions of Gibbs-Markov maps},  Probab. Theory Related Fields   \textbf{123},  no. 1  (2002), 38-40.


\bibitem{Al} J. F. Alves, {\em Strong statistical stability of non-uniformly expanding maps},
Nonlinearity \textbf{17} (2004), 1193–1215.

\bibitem{ABV} J. F. Alves, C. Bonatti, M. Viana, {\em
SRB measures for partially hyperbolic systems whose
       central direction is mostly expanding}, Invent. Math. \textbf{140} (2000), 351-398.



\bibitem{ALP} J. F. Alves, S. Luzzatto, V. Pinheiro, {\em Markov
structures and decay of correlations for non-uniformly expanding
dynamical systems},  Ann. Inst.   Henri Poincaré  (C) Non Linear Anal. \textbf{22}, n.6 (2005) 817-839.



%




%
%
%
%
%

\bibitem{AP} J. F. Alves, V. Pinheiro, {\em
Topological structure of (partially) hyperbolic sets with positive
volume},  Trans. Amer. Math. Soc. \textbf{360} (2008), 5551-5569.

\bibitem{AP2} J. F. Alves, V. Pinheiro,    \emph{ Slow rates of mixing for dynamical systems with hyperbolic structures},
J. Stat. Phys. \textbf{131}, n.3 (2008) 505-534.


\bibitem{BY2} M. Benedicks, L.-S. Young, {\em
Markov extensions and decay of correlations for certain H\'enon
maps}, Ast\'erisque {\bf 261} (2000), 13-56.

\bibitem{BV} C. Bonatti, M. Viana, {\em SRB measures
for partially hyperbolic systems with mostly contracting central
direction}, Israel J. Math. {\bf 115} (2000), 157-193.

\bibitem{Car}
M.~Carvalho,
\newblock \emph{{S}inai-{R}uelle-{B}owen measures for $n$-dimensional derived from
  {A}nosov diffeomorphisms},
\newblock Ergod. Th. \& Dynam. Sys. {\bf13} (1993) 21-44.

\bibitem{Cas}
A.~A. Castro,
\newblock {\em Backward inducing and exponential decay of correlations for
  partially hyperbolic attractors with mostly contracting central
  direction},
\newblock Israel J. of Math. {\bf 130} (2001), 29–75.

\bibitem{D} D. Dolgopyat, \emph{Limit theorems for partially hyperbolic systems},  Trans. Amer. Math. Soc.  \textbf{356},  no. 4 (2004), 1637--1689.


\bibitem{G} S. Gouëzel, {\em Decay of correlations for nonuniformly expanding systems},
 Bull. Soc. Math. France {\bf 134}, n.1
(2006), 1-31.

\bibitem{HM} M. Holland, I. Melbourne, \emph{Central limit theorems and invariance principles for Lorenz attractors},  J. Lond. Math. Soc. (2)  \textbf{76 } (2007),  no. 2, 345-364.

\bibitem{HP} M. W. Hirsch, Morris W.; C. C. Pugh, \emph{Stable manifolds and hyperbolic sets}, 1970  Global Analysis (Proc. Sympos. Pure Math., Vol. XIV, Berkeley, Calif., 1968)  pp. 133--163 Amer. Math. Soc., Providence, R.I.

\bibitem{Ma}    R. Mañé, \emph{Ergodic theory and differentiable dynamics},  Springer-Verlag, Berlin, 1987.

\bibitem{Mel} I. Melbourne, \emph{Large and moderate deviations for slowly mixing dynamical systems}. Proc. Amer. Math. Soc. \textbf{137} (2009) 1735-1741.



\bibitem{MN1} I. Melbourne, M. Nicol, \emph{Large deviations for nonuniformly hyperbolic systems}, Trans. Amer. Math. Soc. \textbf{360} (2008) 6661-6676.

\bibitem{MN2} I. Melbourne, M. Nicol, \emph{A vector-valued almost sure invariance principle for hyperbolic dynamical systems}, Annals of Probability \textbf{37} (2009) 478-505.

\bibitem{PS} Ya. B. Pesin, Ya. G. Sinai,  \emph{Gibbs measures for partially hyperbolic attractors},  Ergodic Theory Dynam. Systems  \textbf{2}, no. 3-4  (1982),  417-438.


\bibitem{T} A. Tahzibi, \emph{Stably ergodic diffeomorphisms which are not partially hyperbolic},  Israel J. Math.  \textbf{142}  (2004), 315--344.

\bibitem{Y1} L.-S. Young, {\em Statistical properties of dynamical
systems with some hyperbolicity}, Ann.  Math.  {\bf 147} (1998),
585-650.

\bibitem{Y2} L.-S. Young, {\em Recurrence times and rates of mixing},
Israel J. Math. {\bf 110} (1999), 153-188.





\end{thebibliography}
\end{document}